
\magnification\magstep1
\baselineskip = 18pt
\def\n{\noindent}
\overfullrule = 0pt
\def\square{\vcenter{\hrule height1pt\hbox{\vrule width1pt height4pt
\kern4pt \vrule width1pt}\hrule height1pt}}
\def\CC{\hbox{{\rm l}\kern-.5em\hbox{\rm C}}}
\def\QQ{\hbox{{\rm l}\kern-.5em\hbox{\rm Q}}}
\def\PP{\hbox{{\rm l}\kern-.2em\hbox{\rm P}}}
\def\RR{\hbox{{\rm l}\kern-.2em\hbox{\rm R}}}
\def\HH{\hbox{{\rm l}\kern-.2em\hbox{\rm H}}}
\def\ZZ{\hbox{{\rm Z}\kern-.5em\hbox{\rm Z}}}

\centerline{\bf The Arithmetic and Geometry of Elliptic Surfaces}
\medskip

\centerline{Peter F. Stiller}\bigskip

\n {\bf Abstract:}\ We survey some aspects of the theory of elliptic
surfaces and give some results aimed at determining the Picard number
of
such a surface. For the surfaces considered, this will be
equivalent to determining the Mordell-Weil
rank of an elliptic curve defined over a function field in one
variable. An
interesting conjecture concerning Galois actions on the relative
de~Rham
cohomology of these surfaces is discussed.\bigskip

This paper focuses on an important class of algebraic surfaces called
elliptic surfaces. The results while geometric in character are
arithmetic
at heart, and  for that reason we devote a fair portion of our
discussion
to
those definitions and facts that make the arithmetic clear.
Later in the paper, we will explain some recent
results and conjectures.    This is a preliminary version, the
detailed
version will appear elsewhere.

There are a number of natural routes
leading to the definition of the class of elliptic surfaces.
Let $E$ denote a compact connected complex manifold with
$\dim _{\CC} E = 2$.

\n {\bf Theorem 1:} (Siegel)\ The field of meromorphic functions on
$E$ has
transcendence degree $\le 2$ over $\CC$, i.e. the field of
meromorphic
functions is:\medskip

\itemitem{1)} $\CC$ constant functions

\itemitem{2)} a finite separable extension of \ $\CC(x)$

\itemitem{3)} a finite separable extension of \ $\CC(x,y)$.$\qquad
\square$\medskip

\n Case 3) is precisely the set of algebraic surfaces, i.e. those
admitting
an
embedding into $\PP^N_{\CC}$. Case 2) was studied by Kodaira, leading
to a
series of three papers:\medskip

\n Kodaira, K., ``On complex analytic surfaces I, II, III,'' Annals
of
Math. 77 and 78, 1963,
 which expound on elliptic surfaces. Kodaira makes the following
definition:\medskip

\n $\underline{\hskip1.5in}$

\n AMS subject classification:\ 14D05, 14J27.

\n Partially supported by ARO grant DAAL 03-88-K-0019.

\n The detailed version of this paper will be submitted for
publication
elsewhere.

\n {\bf Definition 2:}\ $E$ is {\it elliptic\/} if:\medskip

\itemitem{1)} $\exists$ a smooth curve (read Riemann surface) $X$ and
a
proper holomorphic map $\pi\colon E \to X$ of $E$ onto $X$ such
that

\itemitem{2)} $\pi^{-1}(x)$ (with multiplicity) is a non-singular
curve
$E_x$ of genus one, i.e. a torus, for general $x\in X$. (``General''
means
for all but finitely many $x\in X$.)\medskip

\n {\bf Theorem 3:} (Kodaira)\ Transcendence degree $=1$ (case 2)
above)
implies $E$ is elliptic.$\qquad \square$\medskip

\n There are of course many elliptic surfaces which are algebraic and
so
have transcendence degree $=2$.

 From now on $E$ will denote an elliptic surface. One immediate
question is
to determine the nature of the singular fibers of $\pi\colon \
E\to X$
at the finite set of points \hfil\break $S = \{x_1,\ldots, x_n\}
\subset X$ where
the fiber $\pi^{-1}(x_i) = E_{x_i}$ is something other than a
non-singular
curve (occuring with multiplicity one) of genus one.
In the second of the above mentioned papers of Kodaira, a complete
description of the singular fiber types is given. (We will ignore
multiple
fibers, as these can't occur for the type of elliptic surface defined
below.) \medskip

For our purposes, we will narrow the definition of elliptic surface
as
follows: \medskip

\n {\bf Definition 4:}\ A compact connected complex surface $E$ will
be
called {\it elliptic\/} if:\medskip

\itemitem{1)} $\exists$ a smooth curve $X$ (read Riemann surface) and
a
proper holomorphic map $\pi\colon \ E \to X$ mapping $E$ onto
$X$ such that

\itemitem{2)} $\pi^{-1}(x)$ is a non-singular curve of genus one for
general $x\in X$, and

\itemitem{3)} $\pi\colon \ E\to X$ has a section, i.e. $\exists$ a
holomorphic map ${\cal O}\colon \ X\to E$ s.t. $\pi\circ {\cal O} =
1_X$,

\itemitem{4)} $E/X$ is relatively minimal, i.e. there are no
exceptional
curves of the first kind in the fibers,

\itemitem{5)} $E/X$ is not isotrivial.\medskip

\n A few comments are in order. First condition 3) forces $E$ to be
algebraic and devoid of multiple singular fibers. The section ${\cal
O}\colon
\ X\to E$ furnishes a $K(X)$-rational point on the generic fiber
$E^{\rm
gen}$ viewed as a curve over $K(X)$. Thus the generic fiber $E^{\rm
gen}$
is an elliptic curve
over the field $K(X)$ of meromorphic/rational functions on $X$.
Assumption
5) means that we have a non-trivial variation of complex structure in
the
``good'' fibers. Simply put, this means that the $J$-invariant of the
fibers, which can be viewed as a meromorphic/rational function on
$X$, is
non-constant. We denote this function by ${\cal J}\in K(X)$. The fact
that
${\cal J}$ is non-constant allows us, via the Mordell-Weil theorem,
to
conclude that the  group $E^{\rm gen}(K(X))$ of $K(X)$-rational
points
on the generic fiber, or what is the same, the group of sections of
$\pi\colon \ E\to X$, is a finitely generated abelian group. We
denote its
rank by $r_E$.
Finally, condition 4) in the definition implies that we have blown
down all
the exceptional curves in the fibers of $\pi\colon \ E\to X$. $E$ is
then
the unique minimal compactification of the so-called N\'eron model of
the
elliptic curve $E^{\rm gen}/K(X)$.
We remark that the map $\pi$ and the curve $X$ are essentially
uniquely
determined by the fact that the Jacobian of $X$ must be the Albanese
of
$E$.

This places us in a situation analogous to the common arithmetic
situation
where $E$ is an elliptic curve over a number field $K$, where the
N\'eron
model is an arithmetic surface over the ``curve'' ${\rm Spec}({\cal
O}_K)$,
${\cal O}_K$ being the ring of algebraic integers in $K$, and where
the
fiber
over a point in ${\rm Spec}\ {\cal O}_K$ is the reduction of $E$
modulo a
nonzero prime ideal $\varphi \subset {\cal O}_K$. The finitely
generated
abelian group $E(K)$ describes the solutions in $K$ to the
Diophantine
equation(s) defining $E$.

\n Classical Examples:

$$\matrix{{\rm
Legendre}\hfill&Y^2=X(X-1)(X-\lambda)\hfill&\hbox{avoid
characteristic 2}\hfill\cr
\hbox{(Level 2)}&&\hbox{singular fibers at}\ \lambda
=0,1,\infty\hfill\cr
&&\cr
\hbox{Level 3}\hfill&X^3+Y^3 + 1 = \mu XY\hfill&\hbox{avoid
characteristic
3}\hfill\cr
&&\hbox{singular fibers at}\ \mu^3 = 27\ {\rm or}\ \infty\hfill\cr}$$

\n Further Examples (Elliptic Modular Surfaces):

In a paper entitled ``On elliptic modular surfaces'', which appears
in the
Journal Math. Soc. Japan, Vol.~24, No.~1 (1972), T.~Shioda constructs
an
important class of elliptic surfaces. Given a subgroup of finite
index
$\Gamma\subset SL_2(\ZZ)$ with $\left(\matrix{\hfill -1&\hfill 0\cr
\hfill
0&\hfill -1\cr}\right) \notin \Gamma$, Shioda constructs a family of
elliptic curves $E_\Gamma$ over $X_\Gamma$ the modular curve $\Gamma
\backslash \HH	\cup\ \QQ \cup \{\infty\}$ with the obvious monodromy
representation given by $\Gamma$ and where the lattice of periods of
the
fiber over $x\in X_\Gamma$ is homothetic to $\ZZ\tau +\ZZ$ where
$\tau \in
\HH$, the complex upper half-plane, corresponds to $x\in X_\Gamma$.

This leads us naturally into the world of classical automorphic
forms.
We will allude to this in several other places. For now, we content
ourselves with
recalling one of Shioda's results, namely that the space of cusp
forms of
weight three for $\Gamma$, $S_3(\Gamma)$, is naturally isomorphic to
the
space of holomorphic two forms on $E$, $H^0(E, \Omega^2_E$).

We now turn to the main object of interest in this paper:\medskip

\n {\bf Definition 5:}\ The N\'eron-Severi group $NS(E)$ of $E$ is
defined
to
be the group of divisors of $E$ modulo algebraic equivalence (as
opposed to
rational or linear equivalence):
$$NS(E) = {{\rm divisors}\over \hbox{alg. equiv. to 0}} \subset
H^2(E,\ZZ).$$
We remark that for these surfaces algebraic is the same as
homological
equivalence, so the N\'eron-Severi group sits inside $H^2(E,\ZZ)$,
which
can be shown  to be torsion-free. The Picard
number is defined to be
$$\rho_E = {\rm rank}\ NS(E).$$

By the Lefschetz Theorem on (1,1)-classes one also has:

$$\eqalignno{NS(E) = &H^{1,1}\cap H^2(E,\ZZ) \subset H^2(E, \CC)\cr
\noalign{\hbox{or}}
NS(E) = &\hbox{ the group of topological $\CC$-line bundles}\cr
&\hbox{which admit analytic structure}}$$
PROBLEM: \ Calculate $\rho_E$.\medskip

Because we are interested only in the rank of the N\'eron-Severi
group, it
is reasonable to  tensor with $\QQ$ and work with
$$NS(E) \otimes_{\ZZ}\QQ.$$
Now away from the ``bad'' fibers over $S=\{x_1\ldots x_n\}\subset X$
the
family $E|_{X-S} = \pi^{-1}(X-S)$ is locally differentiably trivial,
so it
is natural to use the Leray spectral sequence
$$E^{p,q}_2 = H^p(X, R^q\pi_*\QQ) \Rightarrow H^{p+q}(E,\QQ)$$
to understand $H^2(E,\QQ)$ in terms of the base $X$ and the fibers
which
are tori. The Leray spectral sequence degenerates at $E_2$ and yields
a
filtration
$$0 \subset F^2_{\QQ} \subset F^1_{\QQ} \subset F^0_{\QQ} =
H^2(E,\QQ)$$
where
$$\eqalignno{F^1_{\QQ} = &{\rm ker}(H^2(E,\QQ) \to H^0(X,
R^2\pi_*\QQ))\cr
\noalign{\hbox{consists of classes which restrict to zero on each
fiber,
and}}
 F^2_{\QQ} = &{\rm im}(H^2(X,\QQ) {\buildrel \pi^*\over
\longrightarrow} H^2(E,\QQ)) = \QQ[E_{x_0}]\cr
\noalign{\hbox{is generated by the cohomology class of a fiber.}}}$$
The filtration quotient is \vskip-12pt

$$F^1_{\QQ}/F^2_{\QQ} \cong  H^1(X, R^1\pi_* \QQ).$$
Now the Hodge decomposition on $H^*(E,\CC)$ induces a Hodge structure
on
the filtration quotient above
$$H^1(X,R^1\pi_* \QQ) = H^{2,0} \oplus H^{1,1} \oplus H^{0,2}$$
where $H^{2,0}$ is all of $H^0(E,\Omega^2_E)$ (because the
restriction of a
holomorphic two form to a fiber is necessarily zero when the fiber is
a
curve).

Now there is a well-known
\medskip

\n {\bf Theorem 6:} (Shioda)\ $\rho_E =r_E + 2 + \sum\limits_{s\in S}
(m_s-1)$
where $m_s$ is the number of irreducible components making up the
fiber
$E_s = \pi^{-1}(s)$. Thus the geometric quantity $\rho_E$ is
essentially
the arithmetic quantity $r_E$, the rank of the Mordell-Weil group of
the
generic fiber $E^{\rm gen}$ treated as an elliptic curve over the
field
$K(X)$.$\qquad \square$ \medskip

In practice the numbers $m_s$ are easy to determine, it is $r_E$ that
is in
general impossible to compute.
What should we focus on?

Let $V^i_{\QQ}$ be the span of the algebraic cycles in $F^i_{\QQ}$ so
$$V^i_{\QQ} = (NS(E) \otimes_{\ZZ} \QQ) \cap F^i_{\QQ}$$
and let $W_{\QQ} = V^1_{\QQ}/V^2_{\QQ} = {(NS(E) \otimes_{\ZZ} \QQ)
\cap
F^1_{\QQ} \over F^2_{\QQ}} \subset H^1(X, R^1\pi_* \QQ)$.\medskip

\n {\bf Theorem 7:}\ $\dim_{\QQ} W_{\QQ} = r_E$.$\qquad
\square$\medskip

\n {\bf Proof:}\ See Stiller [5].

It therefore behooves us to look at $R^1\pi_* \QQ$ or $R^1\pi_*\CC$
which
over $X-S$ is a locally constant sheaf of rank two, so that
$$R^1\pi_*\CC|_{X-S}\otimes _{\CC} {\cal O}_{X-S}$$
is a rank two holomorphic vector bundle on $X-S$. Let's study this
bundle.

To get quickly to the heart of the matter we adopt a naive point of
view.
\medskip

\n Pick a base point $x_0 \in X-S$. Here $S$ will contain the support
of
the ``bad'' fibers and some additional points to be named
later.\medskip

\n In a sufficiently small neighborhood $U$ of $x_0$
$$\eqalign{\pi^{-1}&(U) = E|_U\cr
&\downarrow \pi|_{\pi^{-1}(U)}\cr
&U\hfill}$$
$E|_U$ is a ${\cal C}^\infty$-trivial fiber bundle, i.e.
$$\matrix{\pi^{-1}(U)\hfill &\cong&U\times \pi^{-1}(x_0)\cr
\quad \downarrow \pi|_{\pi^{-1}(U)}&&\downarrow pr_1\cr
\quad U\hfill&=&\qquad \ \ U\qquad \qquad .\cr}$$
Write $E_{x_0}$ for $\pi^{-1}(x_0)$ and choose an oriented basis
$\gamma_1,\gamma_2\in H_1(E_{x_0},\ZZ)\cong \ZZ^2 $ for the homology
of the
fiber  and consider
$$\omega_i(x) = \int_{\gamma_i} \Omega|_{E_x}$$
where $\Omega$ is an appropriate meromorphic 1-form on $E$ with poles
only
on the vertical fibers, i.e. $\Omega|_{E^{\rm gen}}$ is a
$K(X)$-rational
differential of the $1^{\rm st}$ kind on the curve $E^{\rm
gen}/K(X)$.
Notice that  for an appropriate finite set of points $S$ which
includes the
support of the singular fibers, the functions $\omega_i(x)$ can be
analytically continued as holomorphic non-vanishing functions
throughout
$X-S$. Moreover,
$${\rm Im}\ \omega_1(x)/\omega_2(x) >0.$$
For a path $\gamma\in \pi_1(X-S, x_0)$, analytic continuation of the
pair
$\omega_1$, $\omega_2$ around $\gamma$ yields
$${\omega_1\choose \omega_2} \longmapsto M_\gamma {\omega_1\choose
\omega_2}$$
where $M_\gamma \in SL_2(\ZZ)$. This is called the {\it monodromy
representation\/} of $E/X$ and the image of the fundamental group in
$SL_2(\ZZ)$ is a subgroup $\Gamma$ of finite index in $SL_2(\ZZ)$
which is
unique up to conjugation in $SL_2(\ZZ)$. This group does not depend
on the
choice of $\Omega$ or on $S$ (provided $S$ contains the support of
the
singular fibers). $\Gamma$ is called the {\it monodromy group\/} of
$SL_2(\ZZ)$.\medskip

\n Recall (see Deligne SLN 163) that the following notions are
equivalent:\medskip

\item{1)} a representation of $\pi_1(X-S, x_0) \longrightarrow
GL_2(\CC)$

\item{2)} a local system (locally constant sheaf) $V$ of rank 2 on
$X-S$

\item{3)} a rank 2 holomorphic vector bundle ${\cal E}_0 =
V\otimes_{\CC}
{\cal O}_{X-S}$ over $X-S$ with integrable holomorphic connection
$D_0$
$${\cal E}_0 {\buildrel D_0\over \longrightarrow} {\cal E}_0 \otimes
_{{\cal O}_{X-S}} \Omega^1_{X-S}$$
having regular singular points.\medskip

\n We recall that  ${\cal E}_0$ can be uniquely extended to a
holomorphic
(algebraic)
rank~2 bundle ${\cal E}$ on $X$ together with a meromorphic
(rational)
connection $D$ having regular singular points. This is known as the
{\it Gauss-Manin connection\/}.\medskip

\item{4)} A second order linear differential operator $\Lambda$
rational/$K(X)$ with regular singular points. In our case $\Lambda
\omega_i = 0$ \quad for $i=1,2$, i.e. $\Lambda$ annihilates the
periods of
$\Omega$ as functions on the base. This is the {\it Picard-Fuchs
equation\/} of
$E/X$ and $\omega_1,\omega_2$ form a basis for the two dimensional
space of
solutions at $x_0$, and elsewhere via analytic continuation.\medskip

\n Our point of view now shifts to this differential equation
($R^1\pi_*\CC|_{X-S}$ is the $V$ above).  We ask ``How much
information
can we recover from $\Lambda$?''.\medskip

\n {\bf Theorem 8:}\ $E$ is determined by $\Lambda$ up to generic
isogeny.
 We remark that $\Lambda$ depends on the choice of $\Omega$ but any
other
choice is $g\Omega$ for $g\in K(X)$ and this transforms $\Lambda$ in
the
obvious simple way. (See Stiller [2]).$\qquad \square$\medskip

\n {\bf Theorem 9:}\ If $E/X$ and $E'/X$ are generically isogeneous
elliptic surfaces over a fixed base curve $X$ then\medskip

\itemitem{1)} $[PSL_2(\ZZ)\colon \overline{\Gamma}] =
[PSL_2(\ZZ)\colon \
\overline{\Gamma}']$

\itemitem{2)} $b_i(E) = b_i(E')$\quad $i=0,\ldots, 4$\ \ \  Betti
numbers

\itemitem{3)} $p_g(E) = p_g(E')$ and $q(E) = q(E') = {\rm genus}\ X$

\itemitem{4)} $\rho_E = \rho_{E'}$

\itemitem{5)} $r_E = r_{E'}$\medskip

\n Remark on the proof:\ 5) is immediate and is used to prove 4) via
the
formula
$$\rho_E = r_E + 2 + \sum_{s\in S}(m_s-1).$$
What is interesting is that $m_s$ is not preserved by generic isogeny
--
only the sum is. In particular a generic isogeny
$$\matrix{E&{\buildrel \phi \over \longrightarrow}&E'\cr
\pi\searrow&&\swarrow \pi'\cr
&X\cr}$$
as a rational map which is an isogeny of the fibers almost
everywhere, may
not extend to all of $E$ as a regular map, and the same may hold for
the
dual isogeny
$$\matrix{E'&{\buildrel \phi'\over \longrightarrow}&E\cr
\pi'\searrow&&\swarrow\pi\cr
&X\cr}\quad .$$
The sum of the $m_s$ can be determined using the exponents at $s$ of
$\Lambda$ (which is easily obtainable local information) independent
of the
isogeny class. For details see Stiller [3].$\qquad \square$\medskip

 Now how do we capture $H^1(X, R^1\pi_*\CC)$ in terms of $\Lambda$?
We
take our clue from Manin. Given any section $s\colon \ X\to E$,
$\pi\circ s
= 1_X$, we can locally (say near $x_0 \in X-S$) take a family of
paths
$\gamma_x$ between the points $s(x)$ and ${\cal O}(x)$ on the fiber
$E_x$
and compute
$$f(x) = \int_{\gamma_x}\Omega|_{E_x}$$
which is defined up to the periods
$$f + m\omega_1 + n\omega_2.$$
Since $\Lambda$ annihilates the periods,
$$\Lambda f = Z$$
turns out to be a well-defined rational function. $f$ is thus
annihilated
by a 3rd order operator $\tilde \Lambda$ and the rank~3 local system
$V_{\tilde \Lambda}$ arises an element of ${\rm
Ext}^1(\underline{\CC},
V_\Lambda)$ where $\underline{\CC}$ is the trivial local system. The
monodromy representation for $\tilde \Lambda$ takes the form
$$\gamma\in \pi_1(X-S,x_0) \longmapsto \left(\matrix{1&m_\gamma\
n_\gamma\cr 0&M_\gamma\cr}\right) \in SL_3(\ZZ).$$
The motivation here is that the cohomology class of the algebraic
cycle
$s-{\cal O}$ is essentially zero on the fibers, and so provides an
element in the Leray filtration quotient
$$F^1_{\CC}/F^2_{\CC } = H^1(X, R^1\pi_*\CC).$$
We make the following definitions:\medskip

\n {\bf Definition 10:}\ $Z\in K(X)$ is {\it exact\/} if $\Lambda f =
Z$
has a global
single-valued meromorphic solution. Thus $Z\in \Lambda K(X) \subset
K(X)$.

\n {\bf Definition 11:}\ $Z \in K(X)$ is {\it locally exact\/} if
$\forall
p \in
X$ the equation $\Lambda f=Z$ restricted to a small neighborhood
$U_p$ of
$p$ has a single-valued meromorphic solution. We denote the set of
locally
exact $Z\in K(X)$ by $L^{\rm para}_\Lambda$. (The notation comes from
the
theory of automorphic forms and the notion of parabolic
cohomology).

\n {\bf Definition 12:}\ H$^1_{IDR}$ is defined to be $L^{\rm
para}_\Lambda/\Lambda K(X)$ and is called the {\it inhomogeneous
de~Rham
cohomology\/}.

 Some remarks are in order. First we have slid over the non-intrinsic
nature of $\Lambda$ which depends on the choice of $\Omega$ and on
the
choice of a derivation ${d\over dx}$ on $K(X)$. A more intrinsic
formulation
would treat $Z$ as $Z(dx)^2$ a meromorphic quadratic differential,
i.e. a
meromorphic section of $(\Omega^1_X)^{\otimes 2}$. In any event
$H^1_{IDR}$
is
independent of any choices. Secondly, local exactness can be
formulated as
a residue condition.\medskip

\n {\bf Theorem 13:}\ $H^1_{IDR}$ is canonically isomorphic to
$H^1(X,
R^1\pi_*\CC)$.\medskip

\n The proof is achieved by showing that both groups are naturally
isomorphic to the subgroup of locally split extensions in ${\rm
Ext}^1(\underline{\CC}, V_\Lambda)$.
 Note that
$${\rm Ext}^1(\underline{\CC}, V_\Lambda) \cong H^1(\pi_1(X-S,
x_0), (V_\Lambda)_{x_0}) \cong H^1(X-S, R^1\pi_*\CC|_{X-S})$$
with the
middle group being the usual group cohomology. The locally split
classes
correspond to parabolic cohomology and $H^1(X, R^1\pi_*\CC)$.
 Here $H^1(X, R^1\pi_*\CC)\hookrightarrow H^1(X_0,
R^1\pi_*\CC|_{X_0})$
where $X_0 = X-S$. This last inclusion comes from the exact sequence
of low
order terms
in the Leray spectral sequence for $i\colon \ X_0\hookrightarrow X$
and
the sheaf $R^1\pi_* \CC|_{X_0}$. (See Stiller [4], [5]).$\qquad
\square$\medskip

\n Now that we have identified $H^1_{IDR}$ with $H^1(X, R^1\pi_*\CC)$
what
about the Hodge decomposition on the latter. \medskip

\n {\bf Theorem 14:}\ There are two divisors ${\cal A}_0 < {\cal A}$
on
$X$, easily
computable in terms of the local behavior of $\Lambda$, such that
every
element of $L({\cal A}_0)$ is locally exact but never exact (unless
0) and
such that no locally exact element in $L({\cal A})\cap L^{\rm
para}_\Lambda$ is ever exact (except 0), and such that
$$L({\cal A}_0) \hookrightarrow H^1_{IDR}$$
corresponds to $H^{2,0}$ in $H^1(X, R^1\pi_*\CC$) and
$$L({\cal A})\cap L^{\rm para}_\Lambda \hookrightarrow H^1_{IDR}$$
corresponds to $H^{2,0} \oplus H^{1,1}$.$\qquad \square$\medskip

\n The point of this result is that we now have
$\underline{\hbox{unique
representatives}}$
of the form $\Lambda f = Z$ for elements of $H^{2,0}$ and
$H^{2,0}\oplus
H^{1,1}$ in $H^1(X, R^1\pi_* \CC)$. (See Stiller [5]).\medskip

\n {\bf Application:}

We assume for the moment that $\left(\matrix{\hfill -1&\hfill 0\cr
\hfill
0&\hfill -1}\right)\notin \Gamma \subset SL_2(\ZZ)$ where $\Gamma$ is
the
monodromy group of $E/X$. We then have	a diagram
$$\matrix{E{\buildrel \sim \over \rightarrow}E_\Gamma \times
_{X_\Gamma} X&
\longrightarrow&E_\Gamma\hfill\cr
\!\!\!\!\!\!\! \pi\!\!\searrow\quad
\downarrow&&\downarrow\pi_\Gamma\hfill\cr
\quad \ \ X&{\buildrel \omega\over
\longrightarrow}&X_\Gamma\hfill\cr}$$
where $\omega = \omega_1/\omega_2$ is the so-called {\it period
map\/}. If
we suppose $X$ is Galois over the modular curve $X_\Gamma$ then $G =
{\rm
Gal}(X/X_\Gamma)$ acts on $H^1(X, R^1\pi_*\QQ)$ and preserves Hodge
type in
$H^1(X, R^1\pi_*\CC)$.\medskip

\n {\bf Problem:}\ Let $V$ be an irreducible rational or complex
representation of $G$. What is the multiplicity of $V$ in
$H^1(X,R^1\pi_*\QQ)$ or $H^1(X,R^1\pi_*\CC)$?\medskip

\n We have been able to show in many cases where $G$ is cyclic, i.e.
$K(X)$ is a cyclic extension of the field of modular functions
$K(X_\Gamma)$, that all multiplicities are one.
 This can't be true in general, but we conjecture that it holds when
$\left(\matrix{\hfill
-1&\hfill 0\cr \hfill 0&\hfill -1\cr}\right) \notin \Gamma$ under
suitable
hypotheses on $G$ (modulo the obvious trivial constituents from
$H^0(X_\Gamma, \Omega^2_{X_\Gamma})$ etc.). When the multiplicities
are all
1, we can explicitly decompose the $G$-modules $H^{2,0}$ and $H^{2,0}
\oplus
H^{1,1}$ using our unique representatives. Since Hodge type is
preserved
and the multiplicities are one, $H^{1,1}$ is the sum of those
irreducible
constituents of $H^{2,0} \oplus H^{1,1}$ not in $H^{2,0}$. If in turn
all
the complex irreducible constituents
(say $G$ is abelian -- so
that over $\CC$ all irreducible constituents $V$ are one dimensional,
and
over $\QQ$ we want eigenvalues which are all primitive
$d^{\rm th}$ roots of one) of a given irreducible rational
representation (dimension $\phi(d)$ in the abelian case) lie in the
$H^{1,1}$ part, we get a contribution (of $\phi(d)$ in the abelian
case) to
$\rho_E$. (See Stiller [5] for examples.)

One approach to the multiplicity problem is suggested by
a similar looking multiplicity problem that goes back to
\medskip

\item{} C. Chevalley and A. Weil, ``\"Uber das Verhalten der
Integrale
ersten Gattung bei Automorphismen des Funktionenk\"orpers,'' Abh.
Math.
Sem. Univ. Hamburg 10 (1934), 358-361.\medskip

\item{} A. Weil, ``\"Uber Matrizenringe auf Riemannschen Fl\"achen
und den
Riemann-Rochschen Satz,'' Abh. Math. Sem. Univ. Hamburg 11 (1936)
110-115.\medskip

\n and in modern exposition:\medskip

\item{} J.F. Glazebrook and D.R. Grayson, ``Galois representations on
holomorphic differentials,'' preprint.\medskip

\n The set-up is

$\widetilde X/\CC$ a curve of genus $\tilde g$

$G$ acts faithfully on $\widetilde X$ so  $G\hookrightarrow {\rm
Aut}(\widetilde X)$

$G$ acts on $H^0(\widetilde X, (\Omega^1_{\widetilde X})^{\otimes
q})$.\medskip

\n The results describe $H^0(\widetilde X, (\Omega^1_{\widetilde
X})^{\otimes
q})$ as a representation of $G$ for $q\ge 1$. Namely  given an
irreducible
complex representation $V$ of $G$, a formula is given, in terms of
local
ramification invariants, for the multiplicity of $V$ in
$H^0(\widetilde X,
(\Omega^1_{\widetilde X})^{\otimes q})$.\bigskip

\centerline{\bf Bibliography}

\item{1.} J. Manin, ``Algebraic curves over fields with
differentiation,''
AMS Transactions, 50 (1966).

\item{2.} P. Stiller, ``Differential equations associated with
elliptic
surfaces,'' Journal Math. Soc. Japan, Vol. 32, No. 2 (1981), pp.
203-233.

\item{3.} P. Stiller, ``Monodromy and invariants of elliptic
surfaces,''
Pacific Journal of Math., Vol. 92, No. 2 (1981), pp. 433-452.

\item{4.} P. Stiller, `` Automorphic forms and the Picard number of
an
elliptic surface,'' Aspects of Mathematics, Vol. E5, Vieweg Verlag,
1984.

\item{5.} P. Stiller, ``The Picard numbers of elliptic surfaces with
many
symmetries,'' Pacific Jour. of Math., Vol. 128, No. 1 (1987).
\bigskip\bigskip

\n Department of Mathematics, Texas A\&M University, College Station,
TX \
77843-3368\medskip

\n E-mail:\ stiller@alggeo.tamu.edu

\end